\documentclass[11pt]{amsart}

\usepackage[a4paper,hmargin=3.5cm,vmargin=3.5cm]{geometry}
\usepackage{amsfonts,amssymb,amscd,amstext}
\usepackage{graphicx}
\usepackage[dvips]{epsfig}

\usepackage{fancyhdr}
\usepackage{libertine}

\usepackage{enumerate}
\usepackage{titlesec}
\usepackage{mathrsfs}
\usepackage{tikz}

\def\r{\mathbb{R}}

\titleformat{\section}
{\filcenter\bfseries\large} {\thesection{.}}{0.2cm}{}
\titleformat{\subsection}[runin]
{\bfseries} {\thesubsection{.}}{0.15cm}{}[.]
\titleformat{\subsubsection}[runin]
{\em}{\thesubsubsection{.}}{0.15cm}{}[.]

\usepackage[up,bf]{caption}

\theoremstyle{definition}

\numberwithin{equation}{section}
\numberwithin{figure}{section}

\usepackage{color}

\begin{document}

\fancyhead[CO]{Mean curvature rigidity } 
\fancyhead[CE]{ R. Souam} 
\fancyhead[RO,LE]{\thepage} 

\thispagestyle{empty}

\vspace*{1cm}
\begin{center}
{\bf\LARGE Mean curvature  rigidity of horospheres, hyperspheres and hyperplanes }

\vspace*{0.5cm}

{\large\bf  Rabah Souam}
\end{center}

\footnote[0]{\vspace*{-0.4cm}

\noindent R. Souam

\noindent Institut de Math\'{e}matiques de Jussieu-Paris Rive Gauche,   UMR 7586, B\^{a}timent Sophie Germain,  Case 7012, 75205  Paris Cedex 13, France.

\noindent e-mail: {\tt rabah.souam@imj-prg.fr}
}

\vspace*{1cm}

\begin{quote}
{\small
\noindent {\bf Abstract}\hspace*{0.1cm} We prove that  horospheres,  hyperspheres and hyperplanes  in a hyperbolic space $\Bbb H^n,\,n\geq 3,$ admit no perturbations with compact support which increase their mean curvature.  This is an extension of  the analogous result in the Euclidean spaces, due to M. Gromov,   which states that a hyperplane in a Euclidean space $\r^n$ admits no compactly supported perturbations having mean curvature $\geq 0.$ \vspace*{0.2cm}

\noindent{\bf Keywords}\hspace*{0.1cm}  mean curvature, mean convexity, tangency principle.

\vspace*{0.2cm}

\noindent{\bf Mathematics Subject Classification (2010)}\hspace*{0.1cm}  53A10, 53C40
}
\end{quote}


The content of the present note  is motivated by the following nice result of M. Gromov \cite{gromov} 

\medskip 

\noindent {\bf Theorem 1.}  {\em
A hyperplane in a Euclidean space $\Bbb R^n$ cannot be perturbed on a compact set so that its mean curvature satisfies $H\geq 0.$
}

\medskip

 Gromov \cite{gromov}  showed that   Theorem 1 can be derived from the non existence of $\Bbb Z^n$-invariant metrics with positive scalar curvature on $\Bbb R^n.$ He also  gave  another direct argument to prove it using  a symmetrization process. 

In what follows, we give another proof  of Theorem 1 and an extension to the hyperbolic spaces using a simple argument.  More precisely,  we  prove the following:

\medskip
 \noindent {\bf Theorem 2.}  {\em
Let $M$ denote a horosphere, a hypersphere or a hyperplane in a hyperbolic space $\Bbb H^n, n\geq 3$ and $H_M\geq 0$ its (constant) mean curvature. Let  $\Sigma$ be a connected properly   embedded $\mathcal C^2$-hypersurface  in  $\Bbb H^n$ which coincides with $M$  outside a compact subset of $\Bbb H^n.$  If  the mean curvature of $\Sigma $ is $ \geq H_M,$ then $\Sigma=M.$ 
 }

\medskip

We recall that a hyperplane in $\Bbb H^n$ is a complete totally geodesic hypersurface and a hypersphere is a connected component  of a set equidistant from a hyperplane. 

The proof uses the tangency principle which goes back to E. Hopf. We recall it with some details in what follows. We first fix some notations and conventions.


Let $\Sigma$ be an  embedded $\mathcal C^2$-hypersurface, with a global unit normal $\nu,$ in a smooth complete $n-$dimensional Riemannian manifold $M.$ We denote by $\sigma$ the second fundamental form of $\Sigma,$ which is defined as follows 
$$\sigma_p(u,v)=- <\nabla_u \nu, v>\quad {\rm for }\, p\in \Sigma,\,u,v \in T_p\Sigma.$$
The shape operator  $S$ of $\Sigma$ is defined as follows
  $$S_pu=-\nabla_u\nu \quad 
{\rm for } \, p\in \Sigma,\, u\in T_p\Sigma,$$  and the (normalized) mean curvature  of $\Sigma$ is the function 
 $$H=\frac{1}{n-1} {\rm trace}\, S.$$
Note that these definitions depend on the choice of the unit normal field $\nu.$ The mean curvature vector field $\vec H:= H\nu$ is, instead, independent of the choice of the unit normal field.  With our conventions, the mean curvature of a unit sphere in the Euclidean space with respect to its interior unit normal is equal to 1. 
\medskip

Let $p\in \Sigma$ and consider  local coordinates $(x_1,\ldots,x_n)$ around $p$ in $M$ so that $T_p\Sigma=\Bbb R^{n-1}\times \{0\}$ and $\frac{\partial}{\partial x_n}(0)=\nu(p).$  An open neighborhood $\mathcal U$ of  $p$ in  $\Sigma$ is  the graph, in the these coordinates, of a $\mathcal C^2$-function $u$ defined on an open neighborhood $\Omega$ of the origin in $\Bbb R^{n-1}\times\{0\}.$ The mean curvature of $\mathcal U,$ computed with  respect to the unit normal field $\nu,$ is given by a nonlinear elliptic operator $$\mathcal M(u)= \mathcal F(x,u,Du,D^2u)$$ 
with $\mathcal F:(x,z,\bf{\xi},\bf{s}\rm) \in \Omega\times\mathcal V \to \mathcal F (x,z,{\bf{\xi}},{\bf {s}})$   a smooth function defined in $\Omega\times \mathcal V$ where $\mathcal V$ is an open subset of $\Bbb R\times\Bbb R^n\times \Bbb R^{n\times n}.$

\medskip

{\bf Definition.} Let  $\Sigma_1$ and $ \Sigma_2$ be  two embedded hypersurfaces in the Riemannian manifold $M$ oriented by global unit normals $\nu_1$ and  $\nu_2$ respectively and let $p$  an   interior point of both $\Sigma_1$ and $\Sigma_2.$  We say that $\Sigma_1\geq \Sigma_2$ near $p$ if the following conditions are satisfied:

(i) $ \nu_1(p)=\nu_2(p),$

(ii) there are local coordinates $(x_1,\dots,x_n)$ around $p$ in $M$ with  $T_p\Sigma_1=T_p\Sigma_2=\Bbb R^{n-1}\times \{0\}$ and $\frac{\partial}{\partial x_n}(0)=\nu_1(p)=\nu_2(p),$  in which $\Sigma_1$ and $\Sigma_2$ are graphs over an open  domain in $\Bbb R^{n-1}\times\{0\}$ of functions $u_1$ and $u_2$ respectively, satisfying $u_1\geq u_2.$

\medskip 
 
Suppose that  $\Sigma_1\geq\Sigma_2$  near $p$ and that their mean curvature functions  (computed with respect to the given normals), in local coordinates as above, verify $\mathcal  M(u_1)\leq \mathcal M(u_2).$ The tangency principle asserts that $\Sigma_1$ coincides with $\Sigma_2$ in a neighborhood of $p.$ The argument goes as follows. 
Consider the difference function $u=u_2-u_1$ and set $u_t=u_1+t(u_2-u_1)$ for $t\in [0,1].$ We can write the following :

\begin{align*} \mathcal M (u_2)-\mathcal M (u_1) = &\left( \int_0^1 \frac{\partial \mathcal F}{\partial z} (u_t) dt\right) u+ \left( \int_0^1 \frac{\partial \mathcal F}{\partial \xi} (u_t) dt\right) \ldotp Du \\
&+\left( \int_0^1 \frac{\partial \mathcal F}{\partial {\bf s\rm} } (u_t) dt\right) \ldotp D^2u.
\end{align*}

Here $\left( \int_0^1 \frac{\partial \mathcal F}{\partial \xi} (u_t) dt\right) \ldotp Du= \sum_i b_i(x) \,u_{x_i} (x)$ with $b_i(x):= \int_0^1 \frac{\partial \mathcal F}{\partial {\xi}_i} (u_t) dt,$ and
$\left( \int_0^1 \frac{\partial \mathcal F}{\partial {\bf s\rm} } (u_t) dt\right) \ldotp D^2u= \sum_{i,j} a_{ij}(x) u_{x_i x_j}(x) $ with $a_{ij}(x):= \int_0^1 \frac{\partial \mathcal F}{\partial {\bf s\rm}_{i,j} } (u_t) dt.$ 

This shows that the $\mathcal C^2$-function  $u=u_2-u_1$ satisfies the   inequality $Lu\geq 0$ where $L$ is the   linear elliptic operator with continuous coefficients defined by 
$$ Lv= \left( \int_0^1 \frac{\partial \mathcal F}{\partial z} (u_t) dt\right) v+ \left( \int_0^1 \frac{\partial \mathcal F}{\partial \xi} (u_t) dt\right) \ldotp Dv +\left( \int_0^1 \frac{\partial \mathcal F}{\partial {\bf s\rm} } (u_t) dt\right) \ldotp D^2v.$$
 By assumption, $u\leq 0$ and  $u$ vanishes at the origin. The strong maximum principle (see, for instance, Chapter 10, Addendum 2, Corollary 19 in \cite {spivak} ) then shows that $u\equiv 0$  in a neighborhood of the origin,  that is, the hypersurfaces $\Sigma_1$ and $\Sigma_2$ coincide in a neighborhood of $p.$  We can thus state the tangency  principle in the following form:  

\medskip

\noindent {\bf Tangency Principle.}  {\em  Let $\Sigma_1$ and $\Sigma_2$ be two embedded $\mathcal C^2$-hypersurfaces in a smooth Riemannian manifold $M, $ oriented by global unit normals. Suppose  the mean curvature functions $H_1$ and $H_2$ of respectively $\Sigma_1$ and $\Sigma_2$  verify  $H_2\geq a\geq H_1$  for some real number $a.$ Let $p$ be an interior 
point of both $\Sigma_1$ and $\Sigma_2.$  If $\Sigma_1\geq\Sigma_2$  near $p$, then $\Sigma_1$ coincides with $\Sigma_2$ in a neighborhood of  $p.$
}

\medskip 

We can now prove Theorem 2.
We will work with the halfspace model  of $\Bbb H^n,$  that is, $$\Bbb H^n= \{(x_1,\ldots,x_n)\in \Bbb R^n, \, x_n>0\}$$
 endowed with the metric $$ds^2=\frac{dx_1^2+\dots+dx_n^2}{x_n^2}.$$

Let $P_0$ be a hyperplane in $\Bbb H^n.$  Without loss of generality, we may assume $P_0$ is the hyperplane $\{x_1=0, x_n>0\}.$ For $t\in\Bbb R,$ we let $P_t$ denote the hyperplane $\{x_1=t, x_n>0\}.$

Let  $\Sigma$ be a properly  embedded  $\mathcal C^2$-hypersurface in $\Bbb H^n$ which coincides with $P_0$ outside a compact subset of $\Bbb H^n.$  Suppose its mean curvature verifies  $H_{\Sigma}\geq 0.$ $\Sigma$ 
separates $\Bbb H^n$ into 2 connected components, a mean convex one $V_+$ and a mean concave one $V_-.$ We may assume 
$V_+$ coincides with the domain $ \{x_1> 0, x_n>0\}$ outside a compact set; the other case being similar.  There is a largest $T\geq 0$ verifying  $\Sigma\cap P_T\neq \emptyset.$ We orient $P_T$ by its unit normal field pointing in the direction  $x_1>0.$ The hypersurfaces $\Sigma$ and $P_T$ are both closed subsets of $\Bbb H^n,$ so their intersection $\Sigma\cap P_T$ is closed in $\Sigma.$ Let us show it is also open in $\Sigma.$ Take 
a point $p\in \Sigma \cap P_T.$  Then the normal to $\Sigma$ at $p$ points in the direction $x_1>0$ (since we assumed $V_+$ coincides with the domain $ \{x_1> 0, x_n>0\}$ outside a compact set)  and $\Sigma\leq P_{T}$ near $p.$ 
As $H_{\Sigma}\geq 0,$ the tangency principle shows that $\Sigma$ coincides with $ P_{T}$    in a neighborhood of $p.$ This shows that  $\Sigma\cap P_{T}$ is open in $\Sigma$ and  since it is also closed, we conclude that $\Sigma=P_{T}$ and $T=0.$

\medskip

Consider now the case of  a horosphere which we may suppose, without loss of generality, is the horosphere $\mathcal H_1= \{x_n=1\}.$ 
We also consider the family of horospheres $\mathcal H_t=\{x_n=t\}, \, t>0,$ having the same asymptotic boundary as $\mathcal H_1.$ We orient  the horospheres $\mathcal H_t$ by their upward pointing unit normal. For this choice of orientation their mean curvature is equal to one. 

  Let 
 $\Sigma$ be a properly embedded   $\mathcal C^2$-hypersurface in $\Bbb H^n$ which coincides with $\mathcal H_1$ outside a compact subset of $\Bbb H^n$ and having mean curvature $\geq 1. $
 $\Sigma$ divides $\Bbb H^n$ into 2 components. The mean convex one coincides with the domain $\{x_n>1\}$ outside a compact set since $\Sigma$ asymptotically coincides with $\mathcal H_1.$    There is a  largest $T\geq 1$ such that $\Sigma\cap \mathcal H_T\neq \emptyset.$
  At a point $p\in \Sigma\cap \mathcal H_T,$ the unit normals to $\Sigma$ and the horopshere 
$\mathcal H_T$ coincide and  $\Sigma\leq \mathcal H_T$ near $p.$ As $H_{\Sigma}\geq 1$ ,  the tangency principle shows that $\Sigma$ coincides with $\mathcal H_T$ in a neighborhood of $p.$ So the set $\Sigma\cap \mathcal H_T$ is open in $\Sigma$
and since it is also closed, we conclude that $\Sigma$ coincides with the horosphere $\mathcal H_T$  and that $T=1.$

\medskip

We consider now the case of hyperspheres.  Let us denote,  this time,  the hyperplane $\{x_1=0, x_n>0\}$ by $S_0.$  For $\beta\in (-\frac{\pi}{2},\frac{\pi}{2}),$ we let $S_{\beta}$ be the intersection with $\Bbb H^n$ of the hyperplane $P_{\beta}$ in $\Bbb R^n$ of equation $\cos\beta\, x_1=\sin\beta \, x_n.$   Note that $\beta$ is the oriented angle between $P_{\beta}$ and  $P_0.$  $S_{\beta}$ is equidistant from $S_0,$  the distance between them being equal to $\log \frac{1+|\sin\beta|}{\cos\beta}.$   We orient $S_\beta$ by the normals pointing in the same direction as the vector $-\cos\beta \,\vec e_1+\sin\beta\,\vec e_n,$ where $\vec e_1,\dots ,\vec e_n$ is the canonical basis in $\Bbb R^n. $ For this choice of orientation the mean curvature  of $S_\beta$ 
is given by  $H_\beta= \sin\beta.$

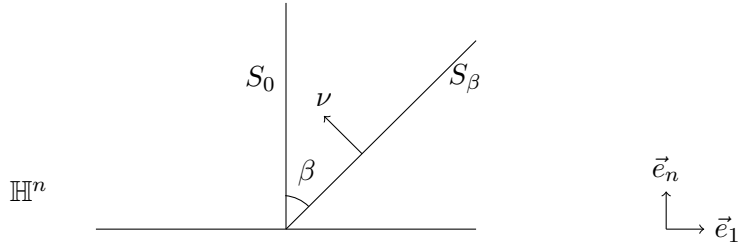
\begin{center}
\begin{tikzpicture}
    \draw (-2.5,0) -- (2.5,0) ;
    \draw (0,0) -- (2.5,2.5);
    \draw (0,0) -- (0,3); 
   \draw (2,2) node[right] {$S_\beta$} ;
    \draw (0.3,0.3) arc (45:85:0.5)  ;
   \draw  (70:0.8) node {$\beta$};
   \draw [->] (1,1) -- (0.5,1.5); 
   \draw (0.5,1.5) node[above] {$\nu$} ;
   \draw (-3,0.5) node[left] {$\Bbb H^n$} ;
    \draw (0,2) node[left] {$S_0$} ;
    
   \draw [->] (5,0) -- (5.5,0) ; 
   \draw [->] (5,0) -- (5,0.5);
   \draw (5,0.5) node[above] {$\vec e_n$} ;
   \draw (5.5,0) node[right] {$\vec e_1$} ;

\end{tikzpicture}
\end{center}

Given a  hypersphere in $\Bbb H^n,$ it is congruent to  $S_{\beta_0}$ for some $\beta_0\in (0, \frac{\pi}{2}).$ Let 
 $\Sigma$ be a properly embedded   $\mathcal C^2$-hypersurface in $\Bbb H^n$ which coincides with $S_{\beta_0} $ outside a compact set and having  mean curvature $H_{\Sigma} \geq \sin\beta_0.$ $\Sigma$ separates $\Bbb H^n$ 
 in 2 components and  the mean convex one  coincides, outside a compact set,  with the mean convex domain bounded by $S_{\beta_0}.$ There is a smallest $\beta\in (-\frac{\pi}{2},\beta_0]$ so that $\Sigma\cap S_\beta \neq \emptyset.$ At a point $p\in \Sigma\cap S_\beta ,$ the unit normals to $S_\beta$ and $\Sigma$ coincide, $\Sigma\leq S_\beta$ near $p$ and 
 $H_\Sigma\geq\sin{\beta_0}\geq \sin\beta=H_\beta.$ By the tangency principle, we know that $\Sigma$ and $S_\beta$ coincide in a neighborhood of $p.$ This shows that the intersection $\Sigma\cap S_\beta$ is open in $\Sigma.$ As it is also closed, we conclude that $\Sigma=S_\beta$ and also that 
 $\beta=\beta_0.$

\medskip
\noindent {\bf Remarks.}  

1. To prove Theorem 1, one uses, as  in the case of a hyperplane in $\Bbb H^n,$ the family of hyperplanes parallel to the given deformed hyperplane.

2. The case of a hyperplane in $\Bbb H^n$ can also be treated using the family of hyperspheres equidistant to it.

3. The argument above can also be utilized to obtain analogous rigidity results for hypersurfaces in other ambient manifolds. Consider, for instance, a product manifold $M\times \Bbb R$ where $M$ is a connected non compact orientable manifold without boundary. 
Suppose $M\times\Bbb R$  is equipped with a Riemannian metric so that $M_t=M\times\{t\}$ has constant mean curvature for each $t\geq 0.$ This happens, for instance, for a warped product metric, that is, a metric of the form $f (t) ds^2_M +dt^2,$ where $ds^2_M$ is a complete metric on $M$ and $f$ a smooth positive function on $[0,+\infty).$ Call  $H_t$ the mean curvature of $M_t$ computed with respect to the unit normal field pointing into the domain $M\times [t,+\infty)$  and suppose  
the function $t\in[0,+\infty) \to H_t$ is non increasing. Then proceeding as above one shows that, for each $t_0\geq 0,$ if $\Sigma$ is a connected  properly embedded  hypersurface which coincides with  $M_{t_0}$  outside a compact subset of $M\times\Bbb R$ and has mean curvature $H_{\Sigma}\geq H_{t_0},$ then $\Sigma=M_{t_0}.$


\end{document}